\documentclass[a4paper,12pt]{amsart}
\usepackage{amsfonts}
\usepackage{amssymb}
\usepackage{ifthen}
\usepackage{graphicx}
\nonstopmode \numberwithin{equation}{section}
\newtheorem{thm}{Theorem}
\newtheorem{cor}{Corollary}
\newtheorem{lem}{Lemma}
\newtheorem{prop}{Proposition}
\newtheorem{opbl}{Open Problem}

\newtheorem{claim}{Claim}
\newtheorem{conj}[equation]{Conjecture}

\theoremstyle{definition}
\newtheorem{defn}{Definition}
\newtheorem{case}{Case}
\newtheorem{examp}[equation]{Example}
\newtheorem{prob}[equation]{Problem}
\newtheorem{ques}[equation]{Question}
\newtheorem{rem}{Remark}

\newcounter {own}
\def\theown {\thesection       .\arabic{own}}

\newenvironment{pf}[1][]{%
 \vskip 3mm
 \noindent
 \ifthenelse{\equal{#1}{}}%
  {{\slshape Proof. }}%
  {{\slshape #1.} }%
 }%
{\qed\bigskip}

\newcounter{alphabet}
\newcounter{tmp}
\newenvironment{Thm}[1][]{\refstepcounter{alphabet}%
\bigskip%
\noindent%
{\bf Theorem \Alph{alphabet}}%
\ifthenelse{\equal{#1}{}}{}{ (#1)}%
{\bf .} \itshape}{\vskip 8pt}

\makeatletter   
\newcommand{\Ref}[1]{\@ifundefined{r@#1}{}{\setcounter{tmp}{\ref{#1}}\Alph{tmp}}}
\makeatother
\newenvironment{Lem}[1][]{\refstepcounter{alphabet}%
\bigskip%
\noindent%
{\bf Lemma \Alph{alphabet}}%
{\bf .} \itshape}{\vskip 8pt}

\newcommand{\IR}{{\mathbb R}}

\newcommand{\diam}{{\operatorname{diam}}}

\def\be{\begin{equation}}
\def\ee{\end{equation}}

\newcommand{\bee}{\begin{enumerate}}
\newcommand{\eee}{\end{enumerate}}

\newcommand{\blem}{\begin{lem}}
\newcommand{\elem}{\end{lem}}
\newcommand{\bthm}{\begin{thm}}
\newcommand{\ethm}{\end{thm}}
\newcommand{\bcor}{\begin{cor}}
\newcommand{\ecor}{\end{cor}}
\newcommand{\beg}{\begin{examp}}
\newcommand{\eeg}{\end{examp}}
\newcommand{\begs}{\begin{examples}}
\newcommand{\eegs}{\end{examples}}
\newcommand{\bdefe}{\begin{defn}}
\newcommand{\edefe}{\end{defn}}
\newcommand{\bprob}{\begin{prob}}
\newcommand{\eprob}{\end{prob}}
\newcommand{\bques}{\begin{ques}}
\newcommand{\eques}{\end{ques}}
\newcommand{\bei}{\begin{itemize}}
\newcommand{\eei}{\end{itemize}}
\newcommand{\bop}{\begin{opbl}}
\newcommand{\eop}{\end{opbl}}
\newcommand{\bcl}{\begin{claim}}
\newcommand{\ecl}{\end{claim}}

\newcommand{\bca}{\begin{case}}
\newcommand{\eca}{\end{case}}
\newcommand{\bcon}{\begin{conj}}
\newcommand{\econ}{\end{conj}}
\newcommand{\bcons}{\begin{conjs}}
\newcommand{\econs}{\end{conjs}}
\newcommand{\bprop}{\begin{prop}}
\newcommand{\eprop}{\end{prop}}
\newcommand{\br}{\begin{rem}}
\newcommand{\er}{\end{rem}}
\newcommand{\brs}{\begin{rems}}
\newcommand{\ers}{\end{rems}}
\newcommand{\bo}{\begin{obser}}
\newcommand{\eo}{\end{obser}}
\newcommand{\bos}{\begin{obsers}}
\newcommand{\eos}{\end{obsers}}
\newcommand{\bpf}{\begin{pf}}
\newcommand{\epf}{\end{pf}}
\newcommand{\ba}{\begin{array}}
\newcommand{\ea}{\end{array}}
\newcommand{\beq}{\begin{eqnarray}}
\newcommand{\beqq}{\begin{eqnarray*}}
\newcommand{\eeq}{\end{eqnarray}}
\newcommand{\eeqq}{\end{eqnarray*}}

\newcommand{\ds}{\displaystyle}

\newcounter{minutes}
\divide\time by 60
\newcounter{hours}
\multiply\time by 60 \addtocounter{minutes}{-\time}

\begin{document}

\bibliographystyle{amsplain}
\title{Local properties of quasihyperbolic and freely quasiconformal
mappings}
\thanks{$^\dagger$ File:~\jobname .tex,
          printed: \number\year-\number\month-\number\day,
          \thehours.\ifnum\theminutes<10{0}\fi\theminutes}
\author{Y. Li}
\address{Y. Li, Department of Mathematics,
Hunan Normal University, Changsha,  Hunan 410081, People's Republic
of China} \email{yaxiangli@163.com}

\author{M. Vuorinen}
\address{M. Vuorinen, Department of Mathematics
University of Turku, 20014 Turku,
Finland
} \email{vuorinen@utu.fi}

\author{X. Wang $^* $
}
\address{X. Wang, Department of Mathematics,
Hunan Normal University, Changsha,  Hunan 410081, People's Republic
of China} \email{xtwang@hunnu.edu.cn}

\date{}
\keywords{ QH mapping, quasisymmetric mapping, quasiconformal mapping, FQC mapping, CQH homeomorphism.\\
${}^{\mathbf{*}}$ Corresponding author}

\begin{abstract}
Suppose that $E$ and $E'$ denote real Banach spaces with dimension at least $2$, that $D\subset E$ and $D'\subset E'$ are
domains, and that $f: D\to D'$ is a homeomorphism. In this paper, we prove that if there exists some constant  $M>1$
(resp. some homeomorphism $\varphi$) such that for all $x\in D$,
$f: B(x,d_D(x))\to f(B(x,d_D(x)))$ is $M$-QH (resp. $\varphi$-FQC), then $f$ is $M_1$-QH with $M_1=M_1(M)$
(resp. $\varphi_1$-FQC with $\varphi_1=\varphi_1(\varphi)$). We apply our results to establish, in terms  of the $j_D$ metric,
a sufficient condition for a homeomorphism to be FQC.

\end{abstract}

\thanks{The research was partly supported by NSF of
China (No. 11071063) and Hunan Provincial Innovation Foundation For Postgraduate and by the Academy of Finland grant of Matti Vuorinen
Project number 2600066611.}

\maketitle\pagestyle{myheadings} \markboth{}{Local properties of quasihyperbolic and freely quasiconformal
mappings}

\section{Introduction and main results}\label{sec-1}

During the past few decades, modern mapping theory and the geometric theory of quasiconformal maps have been studied from several points of view. These studies include Heinonen's work on metric measure spaces \cite{Hei}, Koskela's study of maps with finite distortion \cite{KKM} and V\"ais\"al\"a's work dealing with quasiconformality in infinite dimensional Banach spaces \cite{Vai6-0, Vai6, Vai6', Vai4, Vai8}. The quasihyperbolic metric is an important tool in each of these investigations although their respective methods are otherwise quite divergent. In this paper we will study some questions left open by  V\"ais\"al\"a's work. In passing we remark that the quasihyperbolic geometry has been recently studied by many people(cf. \cite{HIMPS, HPWW, Kle, RT}).

Throughout the paper, we always assume that $E$ and $E'$ denote real Banach
spaces with dimension at least $2$, and that $D\subset E$ and $D'\subset E'$ are
domains. The norm of a vector $z$ in
$E$ is written as $|z|$, and for each pair of points $z_1$, $z_2$ in $E$,
the distance between them is denoted by $|z_1-z_2|$, the closed
line segment with endpoints $z_1$ and $z_2$ by
$[z_1, z_2]$. The
distance from $z\in D$ to the boundary $\partial D$ of $D$ is denoted by $d_D(z)$. For an open ball with center $x$ and radius $r$ we use the notation $\mathbb{B}(x,r)$. The boundary sphere is
denoted by $\mathbb{S}(x,r)$. We begin with the following concepts in line with the notation and terminology
of \cite{TV, Vai2, Vai, Vai6-0, Vai6}.

Let $X$ be a metric space and $\dot{X}=X\cup \{\infty\}$. By a triple in $X$ we mean an ordered sequence
$T=(x,a,b)$ of three distinct points in $X$. The ratio of $T$ is the number $$\rho(T)=\frac{|a-x|}{|b-x|}.$$
If $f: X\to Y$ is  an injective map, the image of a triple  $T=(x,a,b)$  is the triple $fT=(fx,fa,fb)$.

\bdefe \label{def1-0} Let $X$ and $Y$ be two metric spaces,
and let $\eta: [0, \infty)\to [0, \infty)$ be a homeomorphism. Suppose
$A\subset X$. An embedding $f: X\to Y$ is said to be {\it
$\eta$-quasisymmetric} or briefly $\eta$-$QS$, if $\rho(T)\leq t$ implies that $\rho(f(T))\leq \eta(t)$  for each
triple $T$ in $X$ and $t\geq 0$.
We note that $\eta(1)\geq 1$ always holds.\edefe

For convenience, in what follows, we always assume that $x$, $y$, $z$, $\ldots$
denote points in $D$ and $x'$, $y'$, $z'$, $\ldots$
the images in $D'$ of $x$, $y$, $z$, $\ldots$
under $f$, respectively. Also we assume that $\alpha$, $\beta$, $\gamma$, $\ldots$
denote curves in $D$ and $\alpha'$, $\beta'$, $\gamma'$, $\ldots$  the images in $D'$ of
$\alpha$, $\beta$, $\gamma$, $\ldots$
under $f$, respectively.

\bdefe \label{def1-0'} Let  $0<q<1$, and let $D$, $D'$ be metric
spaces in $E$ and $E'$, respectively. A homeomorphism $f:D\to D'$ is
{\it $q$-locally $\eta$-quasisymmetric} if $f|_{\mathbb{B}(a,qr)}$ is
$\eta$-QS whenever $\mathbb{B}(a,r)\subset D$. If $D\not=E$, this
means that $f|_{\mathbb{B}(a, q d_D(a))}$ is $\eta$-QS. When $D=E$, obviously, $f$ is $\eta$-QS. \edefe

\bdefe \label{def1.7'}A map $f:X\to Y$ is uniformly continuous if and only if there is $t_0\in (0,\infty]$ and an embedding $\varphi:[0,t_0)\to [0,\infty)$
with $\varphi(0)=0$ such that $$|f(x)-f(y)|\leq \varphi(|x-y|)$$ whenever $x,y\in X$ and $|x-y|\leq t_0$. We then say that $f$ is $(\varphi, t_0)$-uniformly continuous. If $t_0=\infty$, we briefy say that $f$ is $\varphi$-uniformly continuous.\edefe

The definitions of $k_D$ and $j_D$ metric will be given in section \ref{sec-2}.

\bdefe \label{def1.7} Let $D\not=E$ and $D'\not=E'$ be metric
spaces, and let $\varphi:[0,\infty)\to [0,\infty)$ be a growth
function, that is, a homeomorphism with $\varphi(t)\geq t$. We say
that a homeomorphism $f: D\to D'$ is {\it $\varphi$-semisolid} if
$$k_{D'}(f(x),f(y))\leq \varphi(k_D(x,y))$$
for all $x$, $y\in D$, and {\it $\varphi$-solid} if both $f$ and $f^{-1}$
satisfy this condition.\edefe

The special case $\varphi(t)=Mt,\;M\geq1$, gives the $M$-quasihyperbolic maps or briefly $M$-QH. More precisely, $f$ is called $M$-QH if
$$\frac{k_D(x,y)}{M}\leq k_{D'}(f(x),f(y))\leq Mk_D(x,y)$$ for  all  $x$ and $y$ in $D$.

We say that $f$ is {\it fully $\varphi$-semisolid}
(resp. {\it fully $\varphi$-solid}) if $f$ is
$\varphi$-semisolid (resp. $\varphi$-solid) on every  subdomain of $D$. In particular,
when $D=E$, the subdomains are taken to be proper ones in $D$. Fully $\varphi$-solid mappings are also called {\it freely
$\varphi$-quasiconformal mappings}, or briefly {\it $\varphi$-FQC mappings}.

If $E=\IR^n=E'$, then $f$ is $FQC$ if and only if $f$ is
quasiconformal (cf. \cite{Vai6-0}). See \cite{Vai1, Mvo1} for definitions and
properties of $K$-quasiconformal mappings, or briefly $K$-QC mappings. It is known that each $K$-QC
mapping in
$\IR^n$ is $q$-locally $\eta$-QS for every $q<1$ with $\eta=\eta(K, q, n)$, i.e. $\eta$ depends only on
the constants $K$, $q$ and $n$ (cf. \cite[5.23]{Avv}). Conversely, each
$q$-locally $\eta$-QS mapping in $\IR^n$ is a $K$-QC mapping with
$K=(\eta(1))^{n-1}$ by the metric definition of
quasiconformality (cf. \cite[5.6]{Vai6-0}). Further, in \cite{Vai6-0},
V\"ais\"al\"a proved

\begin{Thm}\label{Thm2'}$($\cite[Theorem 5.10]{Vai6-0}$)$
For a homeomorphism $f: D\to D'$, the following conditions are quantitatively equivalent:

\bee
\item $f$ is $\varphi$-FQC;
\item for some fixed $q\in(0,1)$, both  $f$ and $f^{-1}$  are
locally $\eta-QS$;
\item For every $0<q<1$, there is some  $\eta_q$ such that both $f$ and $f^{-1}$ are  $q$-locally $\eta_q$-QS.
\eee\end{Thm}

For $M$-QH mappings, V\"ais\"al\"a \cite{Vai6-0, Vai8} proved the following.

\begin{Thm}\label{Thm4'}$($\cite[Theorem 4.7]{Vai6-0}$)$ Suppose that $D\not=E$, $D'\not=E'$ and that $f: D\to D'$ is $M$-QH. Then $f$ is fully $4M^2$-QH.\end{Thm}

\begin{Thm}\label{Thm2''}$($\cite[Theorem 5.16]{Vai8}$)$
Suppose that $f: D \to D'$ is a homeomorphism and that each point has a neighborhood $D_1 \subset D$ such that $f_{D_1}: D_1 \to f(D_1)$ is $M$-bilipschitz(or briefly $f$ is locally $M$-bilipschitz). Then $f$ is $M^2$-QH.\end{Thm}

 Recall that $M$-QH need not be bilipschitz. Hence the following problem of V\"ais\"al\"a is natural.

\bop\label{Con1}$($\cite[Section 13]{Vai8}$)$ Suppose that $f: D \to D'$ is a homeomorphism and that
there exists $M >1$ such that for each point
has a neighborhood $D_1 \subset D$ such that $f_{D_1}: D_1 \to f(D_1)$ is $M$-QH. Is $f$ $M'$-QH with $M'=M'(M)$? \eop

The first aim of this paper is to study Open Problem \ref{Con1}. Our result is as follows.

\begin{thm}\label{thm1.0} Suppose that $f: D \to D'$ is a homeomorphism and there exists some constant
$M>1$ such that for each point $x\in D$, $f: \mathbb{B}(x,d_D(x)) \to f(\mathbb{B}(x,d_D(x)))$ is $M$-QH.
Then $f$ is $M_1$-QH with $M_1=M_1(M)$.
\end{thm}

Further, in \cite{Vai6-0},
 V\"ais\"al\"a raised the following open problem.

\bop\label{Con2}$($\cite[Section 7]{Vai6-0}$)$ Suppose that $f: G \to G'$ is a
homeomorphism and that each point has a neighborhood $D \subset G$ such that $f_D: D \to f(D)$ is $\varphi$-FQC.
Is $f$ $\varphi'$-FQC with $\varphi'=\varphi'(\varphi)$? \eop

The second aim of this paper is to consider Open Problem \ref{Con2}.
By applying Theorem \ref{thm1.0}, we will prove the following theorem.

\begin{thm}\label{thm2.0} Suppose that $f: D \to D'$ is a homeomorphism and there exists some homeomorphism
$\varphi$ such that for each point $x\in D$, $f: \mathbb{B}(x,d_D(x)) \to f(\mathbb{B}(x,d_D(x)))$ is
$\varphi$-FQC. Then $f$ is $\varphi_1$-FQC with $\varphi_1=\varphi_1(\varphi)$.
\end{thm}

Applying Theorem \ref{thm2.0} we prove the following.
\begin{thm}\label{thm3.0}Let $\varphi:[0,\infty)\to [0,\infty)$ be homeomorphism with $\varphi(t)\geq t$ for all $t$. Suppose that $f: D \to D'$ is a homeomorphism and that for every subdomain $D_1\subset D,$ we have
\beq\label{eq(00)}\varphi^{-1}(j_{D_1}(x,y))\leq j_{D'_1}(x',y')\leq \varphi(j_{D_1}(x,y))\eeq
with $x,y\in D_1.$ Then $f$ is $\varphi_1$-FQC with $\varphi_1=\varphi_1(\varphi).$ \end{thm}

The following two examples show that the converse of Theorem \ref{thm3.0} is not true.

\begin{examp}\label{exm1}Let $E=R^2$ and $f:D=\mathbb{B}(0,1)\to D'=\mathbb{B}(0,1)\setminus[0,1)$ be a conformal mapping. There exist points $x,y\in D$ such that \eqref{eq(00)} does not hold.
 \end{examp}

\begin{examp}\label{exm2}We consider the broken tube $4.12$ in \cite{Vai6-0}. Let $E$ be an infinite-dimensional separable Hilbert space, and choose an orthonormal base $(e_j)_{j\in \mathbf{Z}}$ of $E$ indexed by the set $\mathbf{Z}$ of all integers. Setting $\gamma'_j=[e_{j-1},e_j]$ we obtain the infinite broken line $\gamma'=\cup\{\gamma'_j:j\in \mathbf{Z}\}$. Let $\gamma$ denote the line spanned by $e_1$, and let $D=\gamma+\mathbb{B}(r)$ with $r\leq \frac{1}{10}$ and $f$ be a locally $M$-bilipschitz homeomorphism from $D$ onto a neighborhood $D'$ of $\gamma'$(For more detail see \cite{Vai6-0}). There exist points $x,y\in D$ such that \eqref{eq(00)} does not hold. \end{examp}

The proofs of Theorem \ref{thm1.0} and Theorems \ref{thm2.0} will be given in Section \ref{sec-3}. The proofs of Theorem \ref{thm3.0} and Examples \ref{exm1} and \ref{exm2} will be given in Section \ref{sec-4}.
 In Section \ref{sec-2}, some necessary preliminaries will be introduced.

\section{Preliminaries}\label{sec-2}

The {\it quasihyperbolic length} of a rectifiable arc or a path
$\alpha$ in the norm metric in $D$ is the number (cf.
\cite{GP,Vai3}):

$$\ell_k(\alpha)=\int_{\alpha}\frac{|dz|}{d_{D}(z)}.
$$

For each pair of points $z_1$, $z_2$ in $D$, the {\it quasihyperbolic distance}
$k_D(z_1,z_2)$ between $z_1$ and $z_2$ is defined in the usual way:
$$k_D(z_1,z_2)=\inf\ell_k(\alpha),
$$
where the infimum is taken over all rectifiable arcs $\alpha$
joining $z_1$ to $z_2$ in $D$. For each pair of points $z_1$, $z_2$ in $D$, we have
(cf. \cite{Vai3})

\beq\label{eq(0000)} \; k_D(z_1,z_2)\geq
\log\Big(1+\frac{|z_1-z_2|}{\min\{d_{D}(z_1), d_{D}(z_2)\}}\Big)\doteq j_D(z_1,z_2).\eeq

Gehring and Palka \cite{GP} introduced the quasihyperbolic metric of
a domain in $\IR^n$. Many of the basic properties of this metric may
be found in \cite{Geo}. Recall that an arc $\alpha$ from $z_1$ to
$z_2$ is a {\it quasihyperbolic geodesic} if
$\ell_k(\alpha)=k_D(z_1,z_2)$. Each subarc of a quasihyperbolic
geodesic is obviously a quasihyperbolic geodesic. It is known that a
quasihyperbolic geodesic between every pair of points in $E$ exists if the
dimension of $E$ is finite, see \cite[Lemma 1]{Geo}. This is not
true in arbitrary spaces (cf. \cite[Example 2.9]{Vai4}).
In order to remedy this shortage, V\"ais\"al\"a introduced the following concepts \cite{Vai6}.

\bdefe \label{def1.4}Let $D\neq E$ and $c\geq 1$. An arc $\alpha\subset D$ is a $c$-neargeodesic if and only if $\ell_k(\alpha[x,y])\leq c\;k_D(x,y)$
for all $x, y\in \alpha$.
\edefe

In \cite{Vai4}, V\"ais\"al\"a proved the following property concerning
the existence of neargeodesics in $E$.

\begin{Thm}\label{LemA} $($\cite[Theorem 3.3]{Vai4}$)$
Let $\{z_1,\, z_2\}\subset D$ and $c>1$. Then there is a
$c$-neargeodesic in $D$ joining $z_1$ and $z_2$.
\end{Thm}

\bdefe \label{def1.6} We say that a homeomorphism $f: D\rightarrow
D'$ is {\it $C$-coarsely $M$-quasihyperbolic}, or briefly
$(M,C)$-CQH, in the quasihyperbolic metric if it satisfies
$$\frac{k_D(x,y)-C}{M}\leq k_{D'}(f(x),f(y))\leq M\;k_D(x,y)+C$$
for all $x$, $y\in D$, where $D$ and $D'$ are domains in $E$ and $E'$,
respectively. \edefe

The following result was proved by V\"ais\"al\"a in \cite{Vai6'}.

\begin{Thm}\label{LemB2} {\rm (\cite[Corollary 2.21]{Vai6'})} Suppose that
$D\subset E$ and $D'\subset E'$ are domains and that $f: D\to D'$ is a homeomorphism.
Then the following conditions are quantitatively equivalent:
\bee
\item $f$ is $\varphi$-FQC;
\item $f$ is $(M,C)$-CQH in every subdomain of $D$,\eee
where $\varphi$ and $(M, C)$ depend only on each other.\end{Thm}

\bdefe \label{def1.8} Let $0<t_0\leq 1$ and let
$\theta:[0,t_0)\to [0,\infty)$ be an embedding with $\theta(0)=0$.
 Suppose that $f: D\to D'$ is a homeomorphism
with $D\not=E$ and $D'\not=E'$.
We say that the homeomorphism $f$ is {\it $(\theta,t_0)$-relative} if
$$\frac{|f(x)-f(y)|}{d_{D'}(f(x))}\leq \theta\Big(\frac{|x-y|}{d_D(x)}\Big)$$
for all $x$, $y\in G$ and $|x-y|<t_0d_D(x)$. If $t_0=1$, then we say that $f$ is
{\it $\theta$-relative}.\edefe

Concerning the relations between relative homeomorphisms and solid mappings,  V\"ais\"al\"a
proved the following.

\begin{Thm}\label{Lem1}$($\cite[Corollary 3.8]{Vai6-0}$)$
Suppose that $D\not= E$, $D'\not= E'$  and $f:D\to D'$ is a
homeomorphism. Then the following conditions are quantitatively
equivalent: \bee \item $f$ and $f^{-1}$ are $\theta$-relative;

\item $f$ and $f^{-1}$ are $(\theta, t_0)$-relative;

\item $f$ is $\varphi$-solid.\eee\end{Thm}

The following lemma is useful for the proof of Theorem \ref{thm1.0}.

\begin{lem}\label{lem-1} Let $D\subset E$ be a domain, and let $x\in D$ and $0<s<1$. If $|x-y|\leq sd_D(x)$, then
$$k_{B_x}(x,y)\leq\frac{1}{1-s}\log\Big(1+\frac{|x-y|}{d_D(x)}\Big),$$ where $B_x=\mathbb{B}(x,d_D(x))$.
\end{lem}

\bpf For each $w\in [x,y]$, we have $d_{B_x}(w)\geq d_D(x)-|x-w|$. Then by the Bernoulli inequality

\begin{eqnarray*}k_{B_x}(x,y)&\leq& \int_{[x,y]}\frac{|dw|}{d_{B_x}(w)}\leq \int^{d_D(x)}_{d_D(x)-|x-y|}\frac{dt}{t}
\\ \nonumber&=& \log\Big(1+\frac{|x-y|}{d_D(x)-|x-y|}\Big)\leq
\log\Big(1+\frac{|x-y|}{(1-s)d_D(x)}\Big)\\ \nonumber &\leq&\frac{1}{1-s}\log\Big(1+\frac{|x-y|}{d_D(x)}\Big),\end{eqnarray*}
from which our lemma follows.\epf

We remark that when $E=R^n$, Lemma \ref{lem-1} coincides with Lemma $3.7$  in \cite{Mvo1}.\medskip

\section{The proofs of Theorems \ref{thm1.0} and
\ref{thm2.0}}\label{sec-3}
\subsection{The proof of Theorem 1}
For fixed $x\in D$, we let $D_1=\mathbb{B}(x,d_D(x))$. Then it follows from
Theorem \Ref{Thm2'} that there exist some constant $c\geq2$ such that $f:D_1 \to D_1'$ is $\frac{c(c+1)}{1+c(c+1)}$ -locally
$\eta$-QS with  $\eta(\frac{1}{c})<1$ and $\eta=\eta(c,M)$.
Theorem \Ref{Thm2'} shows that  $f:D_1 \to D_1'$ is also a $(1-\frac{1}{c^4})$-locally $\eta_1$-QS mapping with $\eta_1=\eta_1(M)$.
Furthermore, we infer from Theorem \Ref{Lem1} that $f:D_1 \to D_1'$ is $\theta$-relative with $\theta=\theta(M)$. Let $q=\frac{c(c+1)}{1+c(c+1)}$,
and let $b_1=4M^2(\theta(q)+\eta_1(1)+\frac{\theta(q)}{1-\eta(\frac{1}{c})})$.

For a fixed $x\in D$, let
$w_1\in \mathbb{S}(x,d_{D}(x))\cap \partial D$. We have the following lemma.

\begin{lem}\label{lem-2}$d_{D'_1}(x')\geq \frac{1}{b_1}d_{D'}(x')$.\end{lem}

\bpf
Suppose on the contrary that
\be\label{eq1}d_{D'_1}(x')< \frac{1}{b_1}d_{D'}(x').\ee
Because $f:D_1 \to D_1'$ is  $\theta$-relative, it follows that for all $z\in \mathbb{S}(x, qd_{D_1}(x))$,
$$\frac{|x'-z'|}{d_{D'_1}(x')}\leq \theta(q).$$
Hence \eqref{eq1} shows
\be\label{eq2}\diam\Big(f\big(\mathbb{S}(x, qd_{D_1}(x))\big)\Big)\leq 2\theta(q)d_{D'_1}(x')\leq\frac{2\theta(q)}{b_1}d_{D'}(x').\ee
Let $w'_2\in \mathbb{S}(x',\frac{1}{2}d_{D'}(x'))\cap f([x,w_1])$ such that $[x,w_2)\subset f^{-1}\Big(\mathbb{B}\big(x',\frac{1}{2}d_{D'}(x')\big)\Big)$. We have the following claim.

{\bf Claim.}
$|w_1-w_2|\geq \frac{1}{c^4}d_D(x)$.

We prove this claim also by contradiction. Suppose that
$$|w_1-w_2|<\frac{1}{c^4}d_D(x).$$

Let $w_3$, $w_4,\cdots, w_{m}$ be the  points in $[w_2,x]$ in the
direction from $w_2$ to $x$ such that for each $i\in \{3,\cdots, m\}$
$$|w_i-w_{i-1}|=c|w_{i-1}-w_{i-2}|\,\; \mbox{and}\,\; |w_{m}-x|<c|w_{m}-w_{m-1}|.$$
It is possible that $w_{m}=x$. Obviously, $k>4$.
Since
$$|w_1-w_{m-1}|>|w_{m-2}-w_{m-1}|=\frac{1}{c}|w_{m}-w_{m-1}|$$ and $$|w_{m-1}-x|=|w_{m-1}-w_{m}|+|w_m-x|<(c+1)|w_{m-1}-w_{m}|,$$ we see that

$$|w_1-w_{m-1}|>\frac{1}{c(c+1)}|w_{m-1}-x|.$$
Then $$|w_{m-1}-x|=|w_1-x|-|w_1-w_{m-1}|<|w_1-x|-\frac{1}{c(c+1)}|w_{m-1}-x|,$$
which implies $$|w_{m-1}-x|<\frac{c(c+1)}{1+c(c+1)}|w_1-x|=\frac{c(c+1)}{1+c(c+1)}d_{D_1}(x).$$
Hence
$$w_{m-1}\in \mathbb{B}\Big(x,\frac{c(c+1)}{1+c(c+1)}d_{D_1}(x)\Big),$$
whence, by (\ref{eq2}),

\be\label{eq3}\max\{|w'_{m-1}-w'_{m}|,|w'_{m}-x'|\}< \frac{2\theta(q)}{b_1}d_{D'}(x').\ee

For each $i\in \{3,\cdots, m-1\}$, let
$y_i=\frac{1}{2}(w_{i-1}+w_{i+1})$. Then
$$|w_1-w_{i-1}|\geq\frac{1}{c}|w_{i-1}-w_i|= \frac{1}{2c}|w_{i-1}-y_i|,$$
which implies $$|w_1-w_{i-1}|\geq \frac{1}{2c+1}d_{D_1}(y_i),$$
whence $$w_{i-1}, w_i, w_{i+1}\in \mathbb{B}\Big(y_i, \frac{2c}{2c+1}d_{D_1}(y_i)\Big)\subset\mathbb{B}(x, qd_{D_1}(y_i)).$$
Because $f$ is $q$-locally $\eta$-QS it follows that
$$\frac{|w'_{i-1}-w'_{i}|}{|w'_i-w'_{i+1}|}\leq \eta\Big(\frac{|w_{i-1}-w_{i}|}{|w_i-w_{i+1}|}\Big)\leq\eta\Big(\frac{1}{c}\Big),$$
which together with the choice of $w_2$ and the inequality (\ref{eq3}) imply that

\begin{eqnarray*}\Big(\frac{1}{2}-\frac{2\theta(q)}{b_1}\Big)d_{D'}(x')&\leq&|w'_2-x'|-|x'-w'_{m}|
\leq\sum_{i=3}^{m}|w'_{i}-w'_{i-1}|
\\ \nonumber&\leq&|w'_{m-1}-w'_{m}|\sum_{i=0}^{m-3}\Big(\eta\big(\frac{1}{c}\big)\Big)^i\leq \frac{2\theta(q)}{\big(1-\eta(\frac{1}{c})\big)b_1}d_{D'}(x')\\ &<&\frac{1}{4}d_{D'}(x').
\end{eqnarray*} This obvious contradiction completes the proof of Claim. 

Hence the above Claim 
shows
$$|x-w_2|=|x-w_1|-|w_1-w_2|\leq \Big(1-\frac{1}{c^4}\Big)d_{D_1}(x).$$
Let $y'\in \mathbb{S}(x',d_{D'_1}(x'))\cap D'_1$, and let $w_0 \in f^{-1}([x',y'))\cap \mathbb{S}(x,|x-w_2|)$.
Because $f$ is also $(1-\frac{1}{c^4})$-locally $\eta_1$-QS, we see that
$$\frac{ b_1}{2}\leq\frac{|x'-w_2'|}{|x'-w'_0|}\leq \eta_1(1).$$
This contradiction completes the proof of Lemma \ref{lem-2}.
\epf

Now we are ready to finish the proof of Theorem \ref{thm1.0}.
\medskip

To prove Theorem \ref{thm1.0}, we only need to prove that for $x,y\in D$, the following  
inequalities hold:
\be\label{eq4}\frac{1}{M_1}k_{D}(x,y)\leq k_{D'}(x',y')\leq M_1 k_{D}(x,y).\ee
We divide the discussions into two cases.

\bca \label{ca1}$|x-y|\leq \frac{1}{2}d_{D}(x)$.\eca

By Lemma \ref{lem-1} we have

\beq\label{l-1} k_{D'}(x',y')\leq k_{D'_1}(x',y')\leq M k_{D_1}(x,y)\leq 2M k_{D}(x,y).\eeq

If $|x'-y'|\leq \frac{1}{2}d_{D_1'}(x')$, we know from Theorem \Ref{Thm4'} that $f^{-1}:D'_2=\mathbb{B}(x',d_{D_1'}(x')) \to f^{-1}(D'_2)$
is $4M^2$-QH. Hence again by Lemma \ref{lem-1} and Lemma \ref{lem-2}
we conclude that

\beq\label{l-2} k_{D}(x,y)&\leq& k_{D_2}(x,y)\leq 4M^2 k_{D'_2}(x',y')\\
\nonumber&\leq&8M^2 \log\big(1+\frac{|x'-y'|}{d_{D_1}(x')}\big)\leq 8b_1M^2 k_{D'}(x',y').\eeq

If $|x'-y'|\geq \frac{1}{2}d_{D_1'}(x')$, then Lemma \ref{lem-2} shows
$$k_{D'}(x',y')\geq \log\Big(1+\frac{|x'-y'|}{d_{D'}(x')}\Big)\geq \log\Big(1+\frac{1}{2b_1}\Big),$$
whence

\beq\label{l-3} k_{D}(x,y)\leq \int_{[x,y]}\frac{|dw|}{d(w)}\leq \log 2\leq\frac{\log 2}{\log\big(1+\frac{1}{2b_1}\big)}k_{D'}(x',y').\eeq

The inequalities \eqref{l-1}, \eqref{l-1} and \eqref{l-3} show that in Case 1 (\ref{eq4}) holds with $M_1=8b_1M^2$.

\bca \label{ca2}$|x-y|>\frac{1}{2}d_{D}(x)$.\eca
It suffices to prove the left side inequality in \eqref{eq4} since the proof for the right one is similar.
It follows from Theorem \Ref{LemA} that there exists a
$2$-neargeodesic $\gamma' $ in $D'$ joining $x'$ and $y'$. Let $x=z_{1}$, and let $z_{2}$ be the first intersection point of $\gamma$
with $\mathbb{S}(z_1,\frac{1}{2}d_{D}(z_1))$ in the direction from
$x$ to $y$. We let $z_{3}$ be the first intersection point of
$\gamma$ with $\mathbb{S}(z_2,\frac{1}{2}d_{D}(z_2))$ in the
direction from $z_{2}$ to $y$. By repeating this procedure, we get a
set $\{z_i\}_{i=1}^p$ of points
 in $\gamma$ such that $y$ is contained
in $\overline{\mathbb{B}}(z_{p},\frac{1}{2}d_{D}(z_p))$, but not in
$\overline{\mathbb{B}}(z_{p-1},\frac{1}{2}d_{D}(z_{p-1}))$. Obviously,
$p> 1$.
Hence Case \ref{ca1} yields
\begin{eqnarray*}
k_{D}(x, y)&\leq& \sum_{i=1}^{p-1} k_{D}(z_i, z_{i+1})+k_{D}(z_p, y)
\\ \nonumber
&\leq& 8b_1M^2\big(\sum_{i=1}^{p-1} \;k_{D'}(z'_i,z'_{i+1})+k_{D'}(z_p',
y')\big)
\\ \nonumber
&\leq& 8b_1M^2\ell_k(\gamma'[x',y'])
\\ \nonumber
&\leq&16b_1M^2\,k_{D'}(x',y').
\end{eqnarray*}
Thus the proof of Theorem \ref{thm1.0} is finished. \qed

\subsection{The proof of Theorem 2}
By Theorem \Ref{LemB2}, we only need to prove that $f$ is fully $(M,C)$-CQH.
For every subdomain $D_1$ in $D$, take $x\in D_1$ and assume that $f: D_2=\mathbb{B}(x,d_{D_1}(x))\to D'_2$ is $\varphi$-FQC.
Then Theorem \Ref{LemB2} shows that
$f: D_2\to D'_2$ is fully $(M_1,C_1)$-CQH, where the constants $M_1$, $C_1$ depend only on $\varphi$.
Hence the similar reasoning as in
 the proof of Theorem \ref{thm1.0} implies that $f: D_1\to D'_1$ is $(M,C)$-CQH, where the constants $M$, $C$
 depend only on $\varphi$. By the arbitrariness of the subdomain $D_1$ in $D$, we see that Theorem \ref{thm2.0} holds. \qed
\bigskip
\section{The proofs of Theorems \ref{thm3.0} and Examples \ref{exm1} and \ref{exm2}}\label{sec-4}

To prove Theorem \ref{thm3.0}, the following definition and theorems are needed.
\bdefe \label{def1.3} A domain $D$ in $E$ is called $c$-{\it
uniform} in the norm metric provided there exists a constant $c$
with the property that each pair of points $z_{1},z_{2}$ in $D$ can
be joined by a rectifiable arc $\alpha$ in $ D$ satisfying (cf. \cite{MS, Vai6})

\bee
\item\label{wx-4} $\ds\min_{j=1,2}\ell (\alpha [z_j, z])\leq c\,d_{D}(z)
$ for all $z\in \alpha$, and

\item\label{wx-5} $\ell(\alpha)\leq c\,|z_{1}-z_{2}|$,
\eee

\noindent where $\ell(\alpha)$ denotes the length of $\alpha$,
$\alpha[z_{j},z]$ the part of $\alpha$ between $z_{j}$ and $z$.
\edefe

In \cite{Vai6},  V\"ais\"al\"a characterized uniform domains by the quasihyperbolic metric.

\begin{Thm}\label{thm0.1} {\rm (\cite[Theorem 6.16]{Vai6})}
For a domain $D$ in $E$, the following are quantitatively equivalent: \bee

\item $D$ is a $c$-uniform domain;
\item $k_D(z_1,z_2)\leq c'\;
 \log\left(1+\frac{\ds{|z_1-z_2|}}{\ds\min\{d_{D}(z_1),d_{D}(z_2)\}}\right)$ for every pair of points $z_1$, $z_2\in D$;
\item $k_D(z_1,z_2)\leq c'_1\;
 \log\left(1+\frac{\ds{|z_1-z_2|}}{\ds\min\{d_{D}(z_1),d_{D}(z_2)\}}\right)+d$ for every pair of points $z_1$, $z_2\in D$.\eee
\end{Thm}

 In the case of domains in $ {\mathbb R}^n \,,$ the equivalence
  of items (1) and (3) in Theorem D is due to Gehring and Osgood \cite{Geo} and the equivalence of items (2) and (3) due to Vuorinen \cite{Mvo2}.

\begin{Thm}\label{ThmO} $($\cite[Lemma 6.7]{Vai6}$)$Suppose that $G$ is a $c$-uniform domain and that $x_0\in G$. Then $G_0=G\setminus\{x_0\}$ is $c_0$-uniform
with $c_0=c_0(c)$.\end{Thm}

\begin{Thm}\label{Thmp} $($\cite[Lemma 3.2]{Vai6-0}$)$Let $X$ be $c$-quasiconvex and let $f:X\to Y$ be a map. Then the following conditions are quantitatively equivalent:\bee
\item $f$ is $(\varphi, t_0)$-uniformly continuous;

\item $f$ is $\varphi$-uniformly continuous;

\item $f$ is $\varphi$-uniformly continuous and there are $M\geq 0$ and $C\geq 0$ such that $\varphi(t)\leq Mt+C$ for all $t$.
\eee\end{Thm}

From the proof of Theorem $5.7$ and summary $5.11$ in \cite{Vai6-0}, we can get the following Lemma.

\begin{Lem}\label{LemO}Suppose that $f:D\to D'$ is a homeomorphism. If for every point $x\in D$, $f:D\setminus\{x\}\to D'\setminus\{x'\}$ is $\varphi$-solid, then $f$ is $\psi$-FQC with $\psi=\psi(\varphi).$\end{Lem}

Now we are ready to prove Theorem \ref{thm3.0}.
\subsection{The proof of Theorem 3} By Theorem \ref{thm2.0} and Lemma \Ref{LemO} we know that to prove the theorem we only need to prove for each $a,b\in D$ with $|a-b|\leq d_D(a)$, $f:G = \mathbb{B}(a, d_D(a))\setminus\{b\}\to G'$ is $\psi$-solid
with $\psi=\psi(\varphi).$
On one hand, choose $0<t_0<1$ such that $\varphi(t_0)\leq \log\frac{3}{2}$. Let $x,y\in G$ be points with $k_G(x,y)\leq t_0$. Then \eqref{eq(00)} gives $|x'-y'|\leq \frac{1}{2}d_{G'}(x')$. Let $B_x=\mathbb{B}(x', d_{G'}(x'))$. Then by  Lemma \ref{lem-1} we obtain

\begin{eqnarray*}k_{G'}(x',y')&\leq& k_{B_x}(x',y')\leq 2 \log\big(1+\frac{|x'-y'|}{d_{G'}(x')}\big)
\\ \nonumber &\leq&2 j_{G'}(x',y')\leq2 \varphi(j_G(x,y))\\ \nonumber &\leq&2\varphi(k_G(x,y)).\end{eqnarray*}
Hence Theorem \Ref{Thmp} yields that $f:G\to G'$ is semi-solid.

On the other hand, Theorem \Ref{ThmO}  show that there exists some constant $c>1$ such that $G$ is $c$-uniform. Hence  Theorem \Ref{thm0.1} yields for each $x,y\in G$
$$k_G(x,y)\leq c' j_G(x,y)\leq c'\varphi (j_{G'}(x',y'))\leq c'\varphi(k_{G'}(x',y')),$$ where $c'$ is a constant depending only on $c$.
The proof of Theorem \ref{thm3.0} is complete.\qed

\subsection{The proof of  Example \ref{exm1}}
By \cite[Lemma 3]{Geo}, we know that conformal mapping is $M$-QH mapping for some constant $M\geq1$. Hence, for each $x,y\in D$, we have
$$\frac{k_D(x,y)}{M}\leq k_{D'}(x',y')\leq Mk_D(x,y).$$
Let $x',y'\in D'$ with $x'=(\frac{1}{2},t)$ and $y'=(\frac{1}{2},-t)$. Then
$$k_{D'}(x',y')\geq \log(1+\frac{1}{t}),$$
and
$$j_{D'}(x',y')=\log\big(1+\frac{|x'-y'|}{\min\{d_{D'}(x'),d_{D'}(y')\}}=\log 3.$$
But

\begin{eqnarray*}j_{D}(x,y)&\geq& \frac{1}{2}k_D(x,y)\geq \frac{1}{2M}k_{D'}(x',y')\geq\frac{1}{2M}\log(1+\frac{1}{t})
\rightarrow\infty,\end{eqnarray*}
as $t\rightarrow0.$
\qed
\subsection{The proof of  Example \ref{exm2}}
By Theorem \ref{Thm2''} that $f$ is $M^2$-QH. Let $x,y\in D$ with $x=\sqrt{2}e_1$, $y=m\sqrt{2}e_1$. Then $d_D(x)=r$, $d_D(y)=r$. Because $f$ is locally $M$-bilipschitz, we get $d_{D'}(x')\geq \frac{r}{M}$ and $d_{D'}(y')\geq \frac{r}{M}$.
Hence the fact ``$D'\subset \mathbb{B}(0,2)$" shows that
 $$j_{D'}(x',y')=\log\big(1+\frac{|x'-y'|}{\min\{d_{D'}(x'),d_{D'}(y')\}}\big)\leq \log(1+\frac{4M}{r}),$$
 but $$j_D(x,y)=\log\big(1+\frac{|x-y|}{\min\{d_{D}(x),d_{D}(y)\}}\big)=\log(1+\frac{\sqrt{2}(m-1)}{r})\rightarrow\infty , $$ as $m\rightarrow\infty$.
 Hence \eqref{eq(00)} does not hold.\qed

\subsection{Remark}
Theorem \ref{thm1.0} and Theorem \ref{thm2.0} imply that the condition ``$f$ is $M$-QH (resp. $\varphi$-FQC)" is quantitatively equivalent to the condition ``for every point $x\in D$, $f$ restricted to the maximal ball $\mathbb{B}(x,d_D(x))$ is $M$-QH (resp. $\varphi$-FQC)". Hence, it is natural to ask if this also holds for $\varphi$-solid.

\bigskip

\noindent {\bf Acknowledgements.} This research was finished when the first author was an academic visitor
in Turku University and the first author was supported by the Academy of Finland grant of Matti Vuorinen
Project number 2600066611. She thanks the
Department of Mathematics in Turku University for hospitality.

\end{document}